\documentclass[a4paper,11pt]{article}
\usepackage{indentfirst,latexsym,bm,color}
\usepackage{amsmath,amssymb,amsfonts}

\usepackage[top=1in,left=1in,right=1in]{geometry}

\makeatletter

\@addtoreset{equation}{section} \makeatother

\begin{document}
\newtheorem{dingli}{Theorem}[section]
\newtheorem{dingyi}[dingli]{Definition}
\newtheorem{tuilun}[dingli]{Corollary}
\newtheorem{zhuyi}[dingli]{Remark}
\newtheorem{yinli}[dingli]{Lemma}
\title{ The minimal and maximal symmetries for $J$-contractive projections
\thanks{ This work was supported by NSF of
China (Nos: 11671242, 11571211) and the
Fundamental Research Funds for the Central Universities
(GK201801011).}}
\author{ Yuan Li\thanks{E-mail address:
 liyuan0401@aliyun.com} , \ Xiaomei Cai \ Jiajia Niu, \ Jiaxin Zhang}
\date{} \maketitle\begin{center}
\begin{minipage}{16cm}
{ \small $$    \ \ \   \ School \ of \ Mathematics \ and \
Information \ Science,\ Shaanxi \ Normal \ University, $$ $$ Xi'an,
710062, People's \ Republic \ of\ China. $$ }
\end{minipage}
\end{center}
 \vspace{0.05cm}
\begin{center}
\begin{minipage}{16cm}
\ {\small {\bf Abstract }
In this note, we firstly consider the structures of symmetries $J$ such that a projection $P$ is $J$-contractive.
Then the minimal and maximal elements of the symmetries $J$ with $P^{\ast}JP\leqslant J$(or $JP\geqslant0)$ are given.
Moreover, some formulas between $P_{(2I-P-P^{\ast})^{+}}$ $(P_{(2I-P-P^{\ast})^{-}})$ and $P_{(P+P^{\ast})^-}$ $(P_{(P+P^{\ast})^+})$ are established.
\\}
\endabstract
\end{minipage}\vspace{0.10cm}
\begin{minipage}{16cm}
\end{minipage}
\end{center}
\begin{center} \vspace{0.01cm}
\end{center}

\section{ Introduction}
Let $\mathcal{H}$ and $\mathcal{K}$ be separable complex Hilbert spaces and
   $\mathcal{B(H)}$  ($\mathcal{B(H,K)}$)
be the set of all bounded linear operators on $\mathcal{H}$ (from
$\mathcal{H}$ to $\mathcal{K}$).  For an operator $A\in
\mathcal{B(H)},$ the adjoint of $A$ is denoted by $A^*$ and $A$ is called a self-adjoint operator if $A=A^*.$ We write
$A\geqslant 0$ if $A$ is a positive operator, meaning $\langle
Ax,x\rangle \geqslant0 $
 for all $x\in\mathcal{H}.$ As usual, the operator order (Loewner partial order) relation $A\geqslant B$ between two self-adjoint operators is defined as $A-
B\geqslant0.$  Also, denote by $\mathcal{B({H})}^{+}$ the set of all positive bounded linear operators on $\mathcal{H}.$ If $A \in\mathcal{B({H})}^{+}$
then $A^{\frac{1}{2}} $ denotes the positive square root of $A.$
 Let $A^+$ and $A^-$ be the positive and negative parts of a self-adjoint operator $A,$
that is, $A^+:=\frac{|A|+A}{2}$ and $A^-:=\frac{|A|-A}{2}.$

For an operator $T\in \mathcal{B(H,K)}, N(T),R(T)$ and $\overline{R(T)}$
denote the null
 space, the range of $T,$ and the closure of $R(T),$  respectively.
An operator $J\in \mathcal{B(H)}$ is said to be a symmetry (or self-adjoint involution) if $J=J^*=J^{-1}.$
In this case, $J^+=\frac{I+J}{2}$ and $J^-=\frac{I-J}{2}$ are
mutually annihilating orthogonal projections. If $J$ is a non-scalar symmetry, then
 an indefinite inner product is defined by \[[x,y]:=\langle Jx, y\rangle \qquad(x,y\in \mathcal{H})\] and
 $(\mathcal{H}, J)$ is called a Krein space ([1]).

Let us denote by $\mathcal{B(H)}^{Id}$ the set of all bounded projections (=idempotents) of $\mathcal{B(H)}.$
A projection $P\in \mathcal{B(H)}^{Id}$ is said to be $J$-projection if $P=JP^{\ast}J.$ In particular,
a projection $P$ is called $J-$positive (or $J-$negative) projection if $JP\geqslant0 $ (or $JP\leqslant0$).
And $P\in \mathcal{B(H)}^{Id}$ is said to be $J$-contractive (or expansive) projection
if $P^{\ast}JP\leqslant J$  (or $P^{\ast}JP\geqslant J).$
 Also, $P_{\mathcal{M}}$ denotes the orthogonal projection onto $\mathcal{M},$
 where $\mathcal{M}$ is a closed subspace of $\mathcal{H}.$
Particularly, we use $P_{A}$ to represent the orthogonal projection onto
$\overline{R(A)}$ and
  $P_{A}^{\perp}:=I-P_{A}.$ Furthermore, $A\simeq B$ means that the operators $A$ and $B$ are unitarily
equivalent, that is  $A=UBU^*,$ for some unitary operators $U.$

It is well-known that an operator $P\in \mathcal{B(H)}^{Id}$ can be written as a $2\times2$ operator matrix:  \begin{equation} P=\left(\begin{array}{cc}I&P_1\\0&0\end{array}\right):R(P)\oplus R(P)^{\perp},\end{equation}
where $P_1\in {\mathcal{B}}(R(P)^{\perp}, R(P)).$
In recent years, the
existence of $J$-positive (negative, contractive, expansive) projections and its properties are considered in [6,7,8]. And some geometry and topological properties of projections and decomposition properties of $J-$projections were studied in [3,5,9-13].
In particular, an exposition of operators in Krein spaces can be found in the lecture by T. Ando [1].

In this note, we firstly consider the structures of symmetries $J$ such that a projection $P$ is $J$-contractive.
In particular, the minimal and maximal elements of the symmetries $J$ with $P^{\ast}JP\leqslant J$ (or $JP\geqslant0)$ are given. That is
$$ min\{J:P^{\ast}JP\leqslant J ,\ J=J^{\ast}=J^{-1}\}=2P_{(P+P^{\ast})^{-}}-I+2P_{N(P+P^{\ast})},$$
$$ max\{J: P^{\ast}JP\leqslant J,\ J=J^{\ast}=J^{-1}\}=2P_{(P+P^{\ast})^{-}}-I+2P_{N(P-P^{\ast})} $$ and
$$max\{J: JP\geqslant0,\ J=J^{\ast}=J^{-1}\}=2P_{(P+P^{\ast})^{+}}-I+2P_{N(P+P^{\ast})}.$$
Moreover, some formulas between $P_{(2I-P-P^{\ast})^{+}}$ $(P_{(2I-P-P^{\ast})^{-}})$ and $P_{(P+P^{\ast})^-}$ $(P_{(P+P^{\ast})^+})$ are established.

\section{ Main results }

  It is well known that every operator $T\in \mathcal{B(H,K)}$ has a (unique) polar decomposition
 $T = U(T^*T)^\frac{1}{2},$
where $U$ is a partial isometry
with kernel space $N(T)$ ([4]).
The following lemmas are needed and their proofs can be found in [8].

{\bf Lemma 1.} ([8, Corollary 11]) Let $P\in \mathcal{B(K,H)}$ and $S=\left(\begin{array}{cc}I&P\\P^{\ast}&0\end{array}\right):\mathcal{H}\oplus\mathcal{K}.$ Then
$$P_{S^-}=\left(\begin{array}{cc}\frac{I-T^{-1}}{2}&-T^{-1}P
\\-P^{\ast}T^{-1}& \frac{V(I+T^{-1})V^*}{2}\end{array}\right):\mathcal{H}\oplus\mathcal{K},$$
where $T=(I+4PP^{\ast})^{\frac{1}{2}}$ and $V$ is the unique partial isometry such that
$P^{\ast}=V(PP^{\ast})^{\frac{1}{2}}$ with $R(V)=\overline{R(P^*)}$ and $R(V^*)=\overline{R(P)}.$

{\bf Lemma 2.} ([8, Lemma 13]) Let $P\in \mathcal{B(H)}^{Id}$ have the form (1.1). If $J$ is a symmetry,
then $P$ is a $J$-projection if and only if $J$ has the operator matrix
\begin{equation}J=\left(\begin{array}{cc}J_1(I+P_1P_1^{\ast})^{-\frac{1}{2}}&J_1(I+P_1P_1^{\ast})^{-\frac{1}{2}}P_1\\
P_1^{\ast}(I+P_1P_1^{\ast})^{-\frac{1}{2}}J_1&J_2(I+P_1^{\ast}P_1)^{-\frac{1}{2}}\end{array}\right):R(P)\oplus R(P)^{\perp},\end{equation}
where $J_1$ and $J_2$ are symmetries on the subspaces
$R(P)$ and $R(P)^{\perp}$, respectively, satisfying $J_1P_1+P_1J_2=0.$

{\bf Lemma 3.} ([8, Corollary 14]) Let $P\in \mathcal{B(H)}^{Id}$ have the form (1.1).\ If $J$ is a symmetry, then
$JP\geqslant0$
if and only if \begin{equation}J=\left(\begin{array}{cc}(I+P_1P_1^{\ast})^{-\frac{1}{2}}&(I+P_1P_1^{\ast})^{-\frac{1}{2}}P_1\\
P_1^{\ast}(I+P_1P_1^{\ast})^{-\frac{1}{2}}&J_2(I+P_1^{\ast}P_1)^{-\frac{1}{2}}\end{array}\right):R(P)\oplus R(P)^{\perp},\end{equation}
where $J_2$ is a symmetry on the subspace
$R(P)^{\perp}$  with $P_1=-P_1J_2.$

{\bf Lemma 4.} ([8, Theorem 15]) Let $P\in \mathcal{B(H)}^{Id}.$ Then  \begin{equation}min\{J: JP\geqslant0,\ J=J^{\ast}=J^{-1}\}=2P_{(P+P^{\ast})^{+}}-I,\end{equation}
where the ``min" is in the sense of Loewner partial order.

In what follows, we give a new characterization for symmetry $J$ with $P^{\ast}JP \leqslant J.$

{\bf Lemma 5.} Let $P\in \mathcal{B(H)}^{Id}$ have the form (1.1). If $J$ is a symmetry,
 then $P^{\ast}JP \leqslant J$ if and only if
\begin{equation}J=\left(\begin{array}{cc}J_1(I+P_1P_1^{\ast})^{-\frac{1}{2}}&J_1(I+P_1P_1^{\ast})^{-\frac{1}{2}}P_1\\
P_1^{\ast}(I+P_1P_1^{\ast})^{-\frac{1}{2}}J_1&(I+P_1^{\ast}P_1)^{-\frac{1}{2}}\end{array}\right):R(P)\oplus R(P)^{\perp},\end{equation}
where $J_1\in \mathcal{B}(R(P))$ is a symmetry  with $J_1P_1+P_1=0.$

{\bf Proof.} Necessity. We assume
$$J=\left(\begin{array}{cc}J_{11}&J_{12}\\J_{12}^{\ast}&J_{22}\end{array}\right): R(P)\oplus R(P)^{\perp},$$
where $J_{11}$ and $J_{22}$ are self-adjoint operators. The equation of
 $J=J^{\ast}=J^{-1}$ implies that \begin{equation}\begin{cases}J_{11}^2+J_{12}J_{12}^{\ast}=I\qquad\qquad \textcircled{1}
&\\J_{11}J_{12}+J_{12}J_{22}=0\;\ \qquad \textcircled{2}&\\J_{12}^{\ast}J_{12}+J_{22}^2=I\qquad\qquad \textcircled{3}\end{cases}\end{equation}

On the other hand, a direct calculation shows that the inequality $P^{\ast}JP \leqslant J$ yields $$\left(\begin{array}{cc}J_{11}&J_{11}P_{1}\\P_{1}^{\ast}J_{11}&P_{1}^{\ast}J_{11}P_{1}\end{array}\right) \leqslant \left(\begin{array}{cc}J_{11}&J_{12}\\J_{12}^{\ast}&J_{22}\end{array}\right).$$
Thus
\begin{equation} J_{22}\geqslant P_{1}^{\ast}J_{11}P_{1} \hbox{ } \hbox{ and }\hbox{   }   J_{12}=J_{11}P_1.\end{equation}
Using equations $ \textcircled{1}$ of (2.5) and (2.6), we have $$J_{11}(I+P_1P_1^{\ast})J_{11}=I,$$
which yields that $J_{11}$ is invertible on the subspace $R(P)$ and $J_{11}^2=(I+P_1P_1^{\ast})^{-1}.$

Setting $J_{1}:=J_{11}(J_{11}^{2})^{-\frac{1}{2}},$ we easily verify that  $J_{1}=J_{1}^{\ast}=J_{1}^{-1}$
and \begin{equation}J_{11}=J_1(I+P_1P_1^{\ast})^{-\frac{1}{2}}.\end{equation}
Also,  equations $ \textcircled{3}$ of (2.5) and (2.6) imply
$$J_{22}^2=I-P_1^*J_{11}^2P_1=I-P_1^*(I+P_1P_1^{\ast})^{-1}P_1=
I-(I+P_1^*P_1)^{-1}P_1^{\ast}P_1=(I+P_1^{\ast}P_1)^{-1}.$$
By equations $\textcircled{2}$ of (2.5) and (2.6), it follows that $J_{11}^{2}P_{1}+J_{11}P_{1}J_{22}=0,$
so $$J_{11}P_{1}+P_{1}J_{22}=0,$$ which induces $P_1^{\ast}J_{11}P_{1}=-P_1^{\ast}P_{1}J_{22}.$
Thus inequality $J_{22}-P_1^{\ast}J_{11}P_{1}\geqslant 0$  implies $$(I+P_1^{\ast}P_1)J_{22}\geqslant 0,$$
which yields $$(I+P_1^{\ast}P_1)J_{22}=(I+P_1^{\ast}P_1)^{\frac{1}{2}}J_{22}(I+P_1^{\ast}P_1)^{\frac{1}{2}}.$$
 Then $J_{22}\geqslant 0,$ so \begin{equation}J_{22}=(I+P_1^{\ast}P_1)^{-\frac{1}{2}}.\end{equation}
Moreover, using equations $\textcircled{2}$ of (2.5) and (2.6)-(2.8), we get that
$$(J_1P_1+P_1)(I+P_1^{\ast}P_1)^{-\frac{1}{2}}=J_{11}P_1+P_1J_{22}=0,$$
that is $J_1P_1+P_1=0.$

Sufficiency is clear from a direct calculation.\qquad $\square$

{\bf Lemma 6.} Let $P\in \mathcal{B(H)}^{Id}$ have the form (1.1). Then

(i)\ \ $N(P+P^{\ast})=0\oplus N(P_{1}).$

(ii) $N(P-P^{\ast})=N(P_{1}^{\ast})\oplus N(P_{1}).$

{\bf Proof.} (i) Clearly, \begin{equation}P+P^{\ast}=\left(\begin{array}{cc}2I&P_{1}\\P_{1}^{\ast}&0\end{array}\right):R(P)\oplus R(P)^{\perp}.\end{equation}
Let $x\in R(P), y\in R(P)^{\perp}$ satisfy $(P+P^{\ast})(x\oplus y)=0.$ Then
$$ 2x+P_{1}y=0 \ \ \hbox{and} \ \ \ P_{1}^{\ast}x=0.$$
So $P_{1}^{\ast}P_{1}y=-2P_{1}^{\ast}x=0,$ which implies that $y\in N(P_{1})$ and $x=0.$ That is
$$N(P+P^{\ast})\subseteq 0\oplus N(P_{1}).$$
Another inclusion relation $0\oplus N(P_{1})\subseteq N(P+P^{\ast})$ is obvious.
Thus $$N(P+P^{\ast})=0\oplus N(P_{1}).$$

(ii) Obviously, $$P-P^{\ast}=\left(\begin{array}{cc}0&P_{1}\\-P_{1}^{\ast}&0\end{array}\right):R(P)\oplus R(P)^{\perp}.$$
Let $x\in R(P)$ and $ y\in R(P)^{\perp}.$ Then $(P-P^{\ast})(x\oplus y)=0$ if and only if
$$P_{1}y=0\ \ \ \hbox{and} \ \ \ P_{1}^{\ast}x=0,$$ which yields $x\in N(P_{1}^{\ast})$ and $y\in N(P_{1}).$
Hence $$N(P-P^{\ast})=N(P_{1}^{\ast})\oplus N(P_{1}).$$  \qquad $\square$

The following is our main result.

{\bf Theorem 7.} Let $P\in \mathcal{B(H)}^{Id}$. Then

 (i) $min\{J: P^{\ast}JP\leqslant J,\ J=J^{\ast}=J^{-1}\}=2P_{(P+P^{\ast})^{-}}-I+2P_{N(P+P^{\ast})}.$

(ii) $max\{J: P^{\ast}JP\leqslant J,\ J=J^{\ast}=J^{-1}\}=2P_{(P+P^{\ast})^{-}}-I+2P_{N(P-P^{\ast})}.$\\
where the ``$min$" and ``$max$" are in the sense of Loewner partial order.

{\bf Proof.}  (i) Suppose that $P$ has the matrix form (1.1). Then  by Lemma 1, we get that
$$P_{(P+P^{\ast})^-}=\dfrac{1}{2}\left(\begin{array}{cc}I-T^{-1}&-T^{-1}P_1
\\-P_1^{\ast}T^{-1}&V(I+T^{-1})V^{\ast}\end{array}\right):R(P)\oplus R(P)^{\perp}.$$
where $T=(I+P_1P_1^{\ast})^{\frac{1}{2}}$ and $V$ is the unique partial isometry such that $P_{1}^{\ast}=V|P_{1}^{\ast}|,$  $R(V)=\overline{R(P_1^*)}$ and $R(V^*)=\overline{R(P_1)}.$
Then by Lemma 6, $$2P_{(P+P^{\ast})^{-}}-I+2P_{N(P+P^{\ast})}=
\left(\begin{array}{cc}-T^{-1}&-T^{-1}P_1\\-P_1^{\ast}T^{-1}&V(I+T^{-1})V^{\ast}-I+ 2P_{N(P_{1})}\end{array}\right).$$
Since $VV^{\ast}=P_{P_{1}^{\ast}},$ then $P_{N(P_{1})}=I-P_{P_{1}^{\ast}}=I-VV^{\ast},$ so
$$V(I+T^{-1})V^{\ast}-I+ 2P_{N(P_{1})}=V(T^{-1}-I)V^{\ast}+I.$$

Clearly,  $P_{1}$ and $V$ have the operator matrices forms
\begin{equation}P_{1}=\left(\begin{array}{cc}0&0\\0&P_{11}\end{array}\right):N(P_{1})\oplus N(P_{1})^{\perp}\rightarrow R(P_{1})^{\perp}\oplus \overline{R(P_{1})}\end{equation}
and $$V=\left(\begin{array}{cc}0&0\\0&V_{1}\end{array}\right): R(P_{1})^{\perp}\oplus \overline{R(P_{1})}\rightarrow N(P_{1})\oplus N(P_{1})^{\perp},$$
respectively, where $V_{1}\in {\mathcal{B}}(\overline{R(P_{1})}, N(P_{1})^{\perp})$ is a unitary operator.
Then
$$V(T^{-1}-I)V^{\ast}+I=\left(\begin{array}{cc}I&0\\0&V_{1}(I+P_{11}P_{11}^{\ast})^{-\frac{1}{2}}V_{1}^{\ast}\end{array}\right)
:N(P_{1})\oplus N(P_{1})^{\perp}$$ and \begin{equation}(I+P_1^{\ast}P_1)^{-\frac{1}{2}}=\left(\begin{array}{cc}I&0\\0&(I+P_{11}^{\ast}P_{11})^{-\frac{1}{2}}\end{array}\right):N(P_{1})\oplus N(P_{1})^{\perp}.\end{equation}
Moreover,  $$V_{1}P_{11}P_{11}^{\ast}V_{1}^{\ast}=V_{1}|P_{11}^{\ast}||P_{11}^{\ast}|V_{1}^{\ast}=P_{11}^{\ast}P_{11},$$ which implies
$V_{1}P_{11}P_{11}^{\ast}=P_{11}^{\ast}P_{11}V_{1},$
so $$V_{1}(I+P_{11}P_{11}^{\ast})^{-\frac{1}{2}}=(I+P_{11}^{\ast}P_{11})^{-\frac{1}{2}}V_{1}.$$
Thus
$$V(T^{-1}-I)V^{\ast}+I=(I+P_{1}^{\ast}P_{1})^{-\frac{1}{2}}.$$

Defining  $$J_{0}:=2P_{(P+P^{\ast})^{-}}-I+2P_{N(P+P^{\ast})}=\left(\begin{array}
{cc}-T^{-1}&-T^{-1}P_1\\-P_1^{\ast}T^{-1}&(I+P_{1}^{\ast}P_{1})^{-\frac{1}{2}}\end{array}\right),$$
we know that $J_{0}$ is a symmetry with $P^{\ast}J_{0}P\leqslant J_{0}$ from Lemma 5.

On the other hand, if $P^{\ast}JP\leqslant J,$ then Lemma 5 implies that $J$ has the form (2.4) and
$$J-J_0
=\left(\begin{array}{cc}(J_{1}+I)(I+P_{1}P_{1}^{\ast})^{-\frac{1}{2}}&(J_{1}+I)(I+P_{1}P_{1}^{\ast})^{-\frac{1}{2}}P_{1}
\\P_{1}^{\ast}(I+P_{1}P_{1}^{\ast})^{-\frac{1}{2}}(J_{1}+I)&0\end{array}\right).$$
It is easy to see that the equation $J_{1}P_{1}=-P_{1}$ yields
$$J_{1}P_{1}P_{1}^{\ast}=-P_{1}P_{1}^{\ast}=-P_{1}(-P_{1}^{\ast}J_{1})=P_{1}P_{1}^{\ast} J_{1},$$
which  implies $$(J_{1}+I)(I+P_{1}P_{1}^{\ast})^{-\frac{1}{2}}
=(J_{1}+I)^\frac{1}{2}(I+P_{1}P_{1}^{\ast})^{-\frac{1}{2}}(J_{1}+I)^\frac{1}{2}\geqslant 0,$$
since $J_{1}+I\geqslant 0$ and $(I+P_{1}P_{1}^{\ast})^{-\frac{1}{2}}\geqslant 0.$
Also, we have $$(J_{1}+I)(I+P_{1}P_{1}^{\ast})^{-\frac{1}{2}}P_{1}=(I+P_{1}P_{1}^{\ast})^{-\frac{1}{2}}(J_{1}P_{1}+P_{1})=0$$
and $$P_{1}^{\ast}(I+P_{1}P_{1}^{\ast})^{-\frac{1}{2}}(J_{1}+I)=0.$$
Hence $J-J_0\geqslant 0,$ which induces the desired result.

(ii) Suppose that $P$ has the matrix form (1.1) and $P^{\ast}JP\leqslant J.$
 Then Lemma 5 implies that $J$ has the form (2.4) and  $J_{1}P_{1}+P_{1}=0,$ so $P_{1}^{\ast}J_{1}+P_{1}^{\ast}=0.$ Thus $P_{1}^{\ast}(I+J_{1})=0,$
which means $$R(I+J_{1})\subseteq N(P_{1}^{\ast}).$$ Then we may assume that
$$I+J_{1}=\left(\begin{array}{cc}J_{11}&J_{12}\\0&0\end{array}\right):N(P_{1}^{\ast})\oplus N(P_{1}^{\ast})^{\perp},$$
 so $J_{12}=0$ follows from the fact that $I+J_{1}=I+J_{1}^{*}.$ Thus
$$I+J_{1}=\left(\begin{array}{cc}J_{11}&0\\0&0\end{array}\right),$$
 which yields \begin{equation}J_{1}=\left(\begin{array}{cc}J_{11}-I&0\\0&-I\end{array}\right):N(P_{1}^{\ast})\oplus N(P_{1}^{\ast})^{\perp}.\end{equation}
Therefore, $$\widetilde{J_{1}}:=max\{J_{1}:J_{1}P_{1}+P_{1}=0,\ J_{1}=J_{1}^{\ast}=J_{1}^{-1}\}
=\left(\begin{array}{cc}I&0\\0&-I\end{array}\right)
=2P_{N(P_{1}^{\ast})}-I.$$

Also, suppose that $P_{1}$ has the operator matrix form (2.10).
Then $$P_{1}P_{1}^{*}=\left(\begin{array}{cc}0&0\\0&P_{11}\end{array}\right)
\left(\begin{array}{cc}0&0\\0&P_{11}^{*}\end{array}\right)
=\left(\begin{array}{cc}0&0\\0&P_{11}P_{11}^{*}\end{array}\right):N(P_{1}^{\ast})\oplus N(P_{1}^{\ast})^{\perp}$$
and
\begin{equation}(I+P_1P_1^{\ast})^{-\frac{1}{2}}=\left(\begin{array}{cc}I&0\\
0&(I+P_{11}P_{11}^{\ast})^{-\frac{1}{2}}\end{array}\right)
:N(P_{1}^{\ast})\oplus N(P_{1}^{\ast})^{\perp}.\end{equation}
Thus $$P_{N(P_{1}^{\ast})}(I+P_1P_1^{\ast})^{-\frac{1}{2}}=P_{N(P_{1}^{\ast})}$$  and
$$P_{N(P_{1}^{\ast})}(I+P_1P_1^{\ast})^{-\frac{1}{2}}P_{1}=P_{N(P_{1}^{\ast})}P_{1}=0.$$

We claim that when setting $J_{1}=\widetilde{J_{1}}=diag(I,-I)$ in equation (2.4),
$J$ is the maximal element with $P^{\ast}JP\leqslant J.$
Indeed, we only need to verify that
$$\left(\begin{array}{cc}(\widetilde{J_{1}}-J_{1})(I+P_1P_1^{\ast})^{-\frac{1}{2}}&(\widetilde{J_{1}}-J_{1})(I+P_1P_1^{\ast})^{-\frac{1}{2}}P_1\\
P_1^{\ast}(I+P_1P_1^{\ast})^{-\frac{1}{2}}(\widetilde{J_{1}}-J_{1})&0\end{array}\right)\geqslant 0,$$
 where $J_1$ has the form (2.12).
Combining the equations (2.12) and (2.13), we get that
$$(\widetilde{J_{1}}-J_{1})(I+P_1P_1^{\ast})^{-\frac{1}{2}}=
\left(\begin{array}{cc}2I-J_{11}&0
\\0&0\end{array}\right)\geqslant 0$$ and
$$(\widetilde{J_{1}}-J_{1})(I+P_1P_1^{\ast})^{-\frac{1}{2}}P_1
=(I+P_1P_1^{\ast})^{-\frac{1}{2}}(\widetilde{J_{1}}-J_{1})P_1=(I+P_1P_1^{\ast})^{-\frac{1}{2}}(-P_1+P_1)=0.$$
So the assertion is valid. Then by (i) and Lemma 6,
$$\begin{array}{rl}&max\{J:P^{\ast}JP\leqslant J ,\ J=J^{\ast}=J^{-1}\}
\\=&\left(\begin{array}{cc}(2P_{N(P_{1}^{\ast})}-I)T^{-1}&(2P_{N(P_{1}^{\ast})}-I)T^{-1}P_1\\
P_1^{\ast}T^{-1}(2P_{N(P_{1}^{\ast})}-I)&(I+P_1^{\ast}P_1)^{-\frac{1}{2}}\end{array}\right)
\\=&\left(\begin{array}{cc}-T^{-1}&-T^{-1}P_1\\
-P_1^{\ast}T^{-1}&(I+P_1^{\ast}P_1)^{-\frac{1}{2}}\end{array}\right)+
\left(\begin{array}{cc}2P_{N(P_{1}^{\ast})}T^{-1}&2P_{N(P_{1}^{\ast})}T^{-1}P_{1}
\\2P_{1}^{\ast}T^{-1}P_{N(P_{1}^{\ast})}&0\end{array}\right)
\\=&(2P_{(P+P^{*})^{-}}-I+P_{N(P+P^{\ast})})+(2P_{N(P_{1}^{\ast})}\oplus0)
\\=&2P_{(P+P^{*})^{-}}-I+2(P_{N(P_{1}^{\ast})}+P_{N(P_{1})})
\\=&2P_{(P+P^{*})^{-}}-I+2P_{N(P-P^{\ast})}.\end{array}$$
where $T=(I+P_1P_1^{\ast})^{\frac{1}{2}}. \qquad \square$

In the following, the maximal element of the symmetries $J$ with $JP\geqslant0$ is given.
Note that the minimal element of those symmetries was obtained in [8, Theorem 15].
Also, we consider the minimal and maximal elements of the symmetries $J$ with $JPJ=P^*.$

{\bf Theorem 8.} Let $P\in \mathcal{B(H)}^{Id}$ . Then

(i) $max\{J: JP\geqslant0,\ J=J^{\ast}=J^{-1}\}=2P_{(P+P^{\ast})^{+}}-I+2P_{N(P+P^{\ast})}.$

(ii) The following statements are equivalent:

\quad (a) $max\{J: JPJ=P^{*},\ J=J^{\ast}=J^{-1}\}$ exists.

\quad (b) $P$ is an orthogonal projection.

\quad (c) $max\{J: JPJ=P^{*},\ J=J^{\ast}=J^{-1}\}=I.$

\quad (d) $min\{J: JPJ=P^{*},\ J=J^{\ast}=J^{-1}\}$ exists.

\quad (e) $min\{J: JPJ=P^{*},\ J=J^{\ast}=J^{-1}\}=-I.$\\

{\bf Proof.} The proof of (i) is similar to that of (ii) of Theorem 7. By Lemma 3, We claim that $$max\{J_{2}: P_{1}=-P_{1}J_{2},\ J_{2}=J_{2}^{\ast}=J_{2}^{-1}\}=P_{N(P_{1})}-P_{N(P_1)^{\perp}}.$$

Indeed, since $P_{1}+P_{1}J_{2}=0,$ then $P_{1}(I+J_{2})=0,$  which implies
  $$I+J_{2}=\left(\begin{array}{cc}J_{11}&J_{12}\\0&0\end{array}\right):N(P_{1})\oplus N(P_{1})^{\perp}.$$
It follows from the fact $J_{2}=J_{2}^{\ast}$ that $J_{12}=0,$
so $$I+J_{2}=\left(\begin{array}{cc}J_{11}&0\\0&0\end{array}\right):N(P_{1})\oplus N(P_{1})^{\perp},$$
which yields that $$J_{2}=\left(\begin{array}{cc}J_{11}-I&0\\0&-I\end{array}\right):N(P_{1})\oplus N(P_{1})^{\perp}.$$
Therefore, $$\widetilde{J_{2}}:=max\{J_{2}: P_{1}=-P_{1}J_{2},\ J_{2}=J_{2}^{\ast}=J_{2}^{-1}\}
=P_{N(P_{1})}-P_{N(P_1)^{\perp}}.$$

Let $P_{1}$ be the operator matrix form (2.10). Then
 $$P_{N(P_{1})}(I+P_1^{\ast}P_1)^{-\frac{1}{2}}=P_{N(P_{1})}.$$ follows from equation (2.11).
Similarly, setting $J_{2}:=\widetilde{J_{2}}$ in equation (2.2), we also get that
$J$ is the maximal element with $JP\geqslant 0.$
Thus by the proof of [8, Theorem 15] (Lemma 4), we have
$$\begin{array}{rl}&max\{J: JP\geqslant0,\ J=J^{\ast}=J^{-1}\}
\\=&\left(\begin{array}{cc}(I+P_1P_1^{\ast})^{-\frac{1}{2}}&(I+P_1P_1^{\ast})^{-\frac{1}{2}}P_1\\
P_1^{\ast}(I+P_1P_1^{\ast})^{-\frac{1}{2}}&(P_{N(P_{1})}-P_{N(P_1)^{\perp}})(I+P_1^{\ast}P_1)^{-\frac{1}{2}}\end{array}\right)
\\=&\left(\begin{array}{cc}(I+P_1P_1^{\ast})^{-\frac{1}{2}}&(I+P_1P_1^{\ast})^{-\frac{1}{2}}P_1\\
P_1^{\ast}(I+P_1P_1^{\ast})^{-\frac{1}{2}}&-(I+P_1^{\ast}P_1)^{-\frac{1}{2}}\end{array}\right)+
\left(\begin{array}{cc}0&0\\0&2P_{N(P_{1})}(I+P_1^{\ast}P_1)^{-\frac{1}{2}}\end{array}\right)
\\=&2P_{(P+P^{*})^{+}}-I+(0\oplus2P_{N(P_{1})})
\\=&2P_{(P+P^{*})^{+}}-I+2P_{N(P+P^{\ast})} .\end{array}$$

(ii) $(a)\Rightarrow (b)$. If $J_{0}:=max\{J: JPJ=P^{*},\ J=J^{\ast}=J^{-1}\},$
then by Lemma 2, $$J_{0}=\left(\begin{array}{cc}J_{01}(I+P_1P_1^{\ast})^{-\frac{1}{2}}&J_{01}(I+P_1P_1^{\ast})^{-\frac{1}{2}}P_1\\
P_1^{\ast}(I+P_1P_1^{\ast})^{-\frac{1}{2}}J_{01}&J_{02}(I+P_1^{\ast}P_1)^{-\frac{1}{2}}\end{array}\right):R(P)\oplus R(P)^{\perp},$$
where $J_{01}$ and $J_{02}$ are self-adjoint involutions on the subspaces
$R(P)$ and $R(P)^{\perp}$, respectively, with $J_{01}P_1+P_1J_{02}=0.$
Clearly, by Lemma 2, for $k=1,2,$ $$\widetilde{J_{k}}:=\left(\begin{array}{cc}(-I_{1})^{k}(I+P_1P_1^{\ast})^{-\frac{1}{2}}&(-I_{1})^{k}(I+P_1P_1^{\ast})^{-\frac{1}{2}}P_1\\
P_1^{\ast}(I+P_1P_1^{\ast})^{-\frac{1}{2}}(-I_{1})^{k}&(-I_{2})^{k+1}(I+P_1^{\ast}P_1)^{-\frac{1}{2}}\end{array}\right)$$
are symmetries and satisfy $\widetilde{J_{k}}P\widetilde{J_{k}}=P^{*},$
where $I_{1},I_{2}$ are identity operators on the subspaces of $R(P)$ and $R(P)^{\perp}$, respectively.

It is easy to verify that $J_{0}\geqslant \widetilde{J_{1}}$ implies $J_{02}\geqslant I_{2}.$ Also, $J_{0}\geqslant \widetilde{J_{2}}$ yields $J_{01}\geqslant I_{1}.$
Thus
$J_{01}=I_{1}$ and $ J_{02}=I_{2},$ so $P_{1}=0$ follows from the fact $J_{01}P_1+P_1J_{02}=0.$
Hence, $P$ is an orthogonal projection. In a similarly way, we have $(d)\Rightarrow (b).$

$(b)\Rightarrow (c)\Rightarrow (a)$ and $(b)\Rightarrow (e)\Rightarrow (d)$ are obvious.
\qquad $\square$

{\bf Remark.} Let $P\in \mathcal{B(H)}^{Id}.$ According to the proof of [8, Theorem 15] and [3, Proposition 6.2], we know that
$$2P_{(P+P^{\ast})^{+}}-I=\left(\begin{array}{cc}(I+P_1P_1^{\ast})^{-\frac{1}{2}}&(I+P_1P_1^{\ast})^{-\frac{1}{2}}P_1\\
P_1^{\ast}(I+P_1P_1^{\ast})^{-\frac{1}{2}}&-(I+P_1^{\ast}P_1)^{-\frac{1}{2}}\end{array}\right)
=(P+P^*-I)|P+P^*-I|^{-1}.$$
Thus $$max\{J: JP\geqslant0,\ J=J^{\ast}=J^{-1}\}= (P+P^*-I)|P+P^*-I|^{-1}+2P_{N(P+P^{\ast})}$$
and $$(P+P^*-I)|P+P^*-I|^{-1}P_{N(P+P^{\ast})}=-P_{N(P+P^{\ast})}.$$ \qquad $\square$

In the following, we shall study the relations between $P_{(2I-P-P^{\ast})^{+}}$ $(P_{(2I-P-P^{\ast})^{-}})$ and $P_{(P+P^{\ast})^-}$ $(P_{(P+P^{\ast})^+}).$ The following observation is needed.

{\bf Proposition 9.} Let $P\in \mathcal{B(H)}^{Id}.$ If $P$ and $I-P$ have respectively the operator matrices
\begin{equation}P=\left(\begin{array}{cc}I&P_{1}\\0&0\end{array}\right):R(P)\oplus R(P)^{\perp}\ \hbox{and} \
I-P=\left(\begin{array}{cc}I&Q_{1}\\0&0\end{array}\right):R(I-P)\oplus R(I-P)^{\perp},\end{equation}
then there exist unitary operators $U_{1}\in {\mathcal{B}}(R(P)^{\perp},R(I-P))$ and
$V_{1}\in {\mathcal{B}}(R(I-P)^{\perp},R(P))$ such that $Q_{1}=U_{1}P_{1}^{\ast}V_{1}.$

{\bf Proof.} It is clear that $$I-P=\left(\begin{array}{cc}0&-P_{1}\\0&I\end{array}\right):R(P)\oplus R(P)^{\perp}.$$
Thus there exists a unitary operator $$\widetilde{U}=\left(\begin{array}{cc}U_{11}&U_{12}\\U_{21}&U_{22}\end{array}\right):
R(P)\oplus R(P)^{\perp}\rightarrow R(I-P)\oplus R(I-P)^{\perp},$$
such that $$\left(\begin{array}{cc}U_{11}&U_{12}\\U_{21}&U_{22}\end{array}\right)\left(\begin{array}{cc}0&-P_{1}\\0&I\end{array}\right)
=\left(\begin{array}{cc}I&Q_{1}\\0&0\end{array}\right)\left(\begin{array}{cc}U_{11}&U_{12}\\U_{21}&U_{22}\end{array}\right),$$
which yields $U_{11}+Q_{1}U_{21}=0$ and $-U_{21}P_{1}+U_{22}=0.$ So
\begin{equation}U_{11}=-Q_{1}U_{21}\end{equation} and \begin{equation}U_{22}=U_{21}P_{1}.\end{equation}

On the other hand, the equation $\widetilde{U^{\ast}}\widetilde{U}=I$ implies that
$$U_{11}^{\ast}U_{11}+U_{21}^{\ast}U_{21}=I,$$  so \begin{equation}U_{21}^{\ast}(I+Q_1^{\ast}Q_1)U_{21}=I\end{equation}follows from equation (2.15).

Analogously, the equation $\widetilde{U}\widetilde{U^{\ast}}=I$ implies that
$$U_{21}U_{21}^{\ast}+U_{22}U_{22}^{\ast}=I,$$ so \begin{equation}U_{21}(I+P_1P_1^{\ast})U_{21}^{\ast}=I\end{equation} follows from equation (2.16).

In view of (2.17) and (2.18), $U_{21}$ is invertible from the subspace $R(P)$ onto $R(I-P)^{\perp}.$
Using the polar decomposition theorem, we denote \begin{equation}U_{21}=U(U_{21}^{\ast}U_{21})^{\frac{1}{2}},\end{equation}
where $U$ is a unitary operator from the subspace $R(P)$ onto $R(I-P)^{\perp}.$
Combining equations (2.18) and (2.19), we have
\begin{equation}(U_{21}^{\ast}U_{21})^{\frac{1}{2}}=(I+P_1P_1^{\ast})^{-\frac{1}{2}}.\end{equation}
Also, equation (2.19) implies $U_{21}=(U_{21}U_{21}^{\ast})^{\frac{1}{2}}U,$ so
\begin{equation}(U_{21}U_{21}^{\ast})^{\frac{1}{2}}=(I+Q_1^{\ast}Q_1)^{-\frac{1}{2}}\end{equation}
follows from equation (2.17).

In a similar way, equations $\widetilde{U}\widetilde{U^{\ast}}=I$ and $\widetilde{U^{\ast}}\widetilde{U}=I$  imply
$$U_{11}U_{11}^{\ast}+U_{12}U_{12}^{\ast}=I \hbox{   }\hbox{ and }\hbox{ } U_{12}^{\ast}U_{12}+U_{22}^{\ast}U_{22}=I.$$

Thus equations (2.15) and (2.21) yield that \begin{equation}U_{12}U_{12}^{\ast}=I-Q_1(U_{21}U_{21}^{\ast})Q_{1}^{\ast}
=I-Q_1(I+Q_1^{\ast}Q_1)^{-1}Q_{1}^{\ast}=(I+Q_1Q_1^{\ast})^{-1}.\end{equation}
Similarly, equations (2.16) and (2.21)  imply
\begin{equation}U_{12}^{\ast}U_{12}=I -U_{22}^{\ast}U_{22}=I-P_1^*(I+P_1P_1^{\ast})^{-1}P_{1}=(I+P_1^{\ast}P_1)^{-1}.\end{equation}
Thus $U_{12}$ is invertible from the subspace $R(P)^{\perp}$ onto $R(I-P).$
Using the polar decomposition theorem again, we get that
\begin{equation}U_{12}=V(U_{12}^{\ast}U_{12})^{\frac{1}{2}}=V(I+P_1^{\ast}P_1)^{-\frac{1}{2}},\end{equation}
where $V$ is a unitary operator from the subspace $R(P)^{\perp}$ onto $R(I-P).$

Furthermore, equation $\widetilde{U}\widetilde{U^{\ast}}=I$ implies that
$$U_{21}U_{11}^{\ast}+U_{22}U_{12}^{\ast}=0.$$
Combining equations (2.15), (2.16), (2.19), (2.20) and (2.24), we conclude that
$$-U(I+P_1P_1^{\ast})^{-1}U^{\ast}Q_{1}^{\ast}+U(I+P_1P_1^{\ast})^{-\frac{1}{2}}P_{1}(I+P_1^{\ast}P_1)^{-\frac{1}{2}}V^{\ast}=0.$$
Then $$U(I+P_1P_1^{\ast})^{-1}(-U^{\ast}Q_{1}^{\ast}+P_{1}V^{\ast})=0,$$
so $Q_{1}=VP_{1}^{\ast}U^{\ast}.$
 Setting $U_{1}=V$ and $ V_{1}=U^{\ast},$
  we get that $U_{1}\in {\mathcal{B}}(R(P)^{\perp},R(I-P))$ and
$V_{1}\in {\mathcal{B}}(R(I-P)^{\perp},R(P))$ are unitary operators with $Q_{1}=U_{1}P_{1}^{\ast}V_{1}.$   \qquad $\square$

{\bf Corollary 10.} Let $P\in \mathcal{B(H)}^{Id}.$ Then

(i)\ \ \ $P\simeq P^{\ast}.$

(ii)\ \ $R(P+P^{\ast})$ is closed if and only if $R(2I-P-P^{\ast})$ is closed.

(iii) $P+P^{\ast}+2P_{P}^{\perp}\simeq 2I-P-P^{\ast}+2P_{(I-P)}^{\perp}.$

{\bf Proof.} Let $P$ and $I-P$ have the operator matrices (2.14). Then by Proposition 9, we have $$P_{1}=V_{1}Q_{1}^{\ast}U_{1},$$
where unitary operators $U_{1}\in {\mathcal{B}}(R(P)^{\perp},R(I-P))$ and $V_{1}\in {\mathcal{B}}(R(I-P)^{\perp},R(P)).$

(i) Obviously, $$P^{\ast}=I-(I-P)^*=\left(\begin{array}{cc}0&0\\-Q_{1}^{\ast}&I\end{array}\right):R(I-P)\oplus R(I-P)^{\perp}.$$ Defining a unitary operator $$U:=\left(\begin{array}{cc}0&-U_{1}\\V_{1}^{\ast}&0\end{array}\right):
R(P)\oplus R(P)^{\perp}\rightarrow R(I-P)\oplus R(I-P)^{\perp},$$
 we easily verify that $U^{\ast}P^{\ast}U=P,$ that is $P\simeq P^{\ast}.$

(ii) Setting $$S:=\left(\begin{array}{cc}I&-\frac{P_{1}}{2}\\0&I\end{array}\right):R(P)\oplus R(P)^{\perp},$$
we know that $S$ is invertible.
Considering that $P+P^{\ast}$ has form (2.9) and
$$S^{\ast}(P+P^{\ast})S=\left(\begin{array}{cc}2I&0\\0&-\frac{P_{1}^{\ast}P_{1}}{2}\end{array}\right),$$
 we get that $R(P+P^{\ast})$ is closed if and only if  $R(P_{1}^{\ast}P_{1})$ is closed.
  By a well-known theorem ( or see [14, Lemma 2.1]),
  $R(P_{1}^{\ast}P_{1})$ is closed if and only if $R(P_{1}^*)$ is closed if and only if $R(P_{1})$ is closed.
  Thus $R(P+P^{\ast})$ is closed if and only if $R(P_{1})$ is closed.
Similarly, $R(2I-P-P^{\ast})$ is closed if and only if $R(Q_{1}^{\ast})$ is closed.
Moreover, Proposition 9 implies $R(P_{1})$ is closed if and only if $R(Q_{1}^{\ast})$ is closed,
so the desired result holds.

(iii) It is easy to check that \begin{equation}P+P^{\ast}+2P_{P}^{\perp}=\left(\begin{array}{cc}2I&P_{1}
\\P_{1}^{\ast}&2I\end{array}\right):R(P)\oplus R(P)^{\perp}.\end{equation}
Defining a unitary operator \begin{equation}\widetilde{U}:=\left(\begin{array}{cc}0&V_{1}\\U_{1}^{\ast}&0\end{array}\right):
R(I-P)\oplus R(I-P)^{\perp}\rightarrow R(P)\oplus R(P)^{\perp},\end{equation}
we get that
$$ \widetilde{U}^{\ast}(P+P^{\ast}+2P_{P}^{\perp})\widetilde{U}=
\left(\begin{array}{cc}2I&U_{1}P_{1}^{\ast}V_{1}\\V_{1}^{\ast}P_{1}U_{1}^{\ast}&2I\end{array}\right)
=\left(\begin{array}{cc}2I&Q_{1}\\Q_{1}^{\ast}&2I\end{array}\right)
=2I-P-P^{\ast}+2P_{(I-P)}^{\perp}.$$ \qquad $\square$

{\bf Lemma 11.} ([2, Lemma 2.10] or [6, Proposition 5]) Let $P\in \mathcal{B(H)}^{Id}.$ Then
$P^{\ast}JP\leqslant J$  if and only if $J(I-P)\geqslant 0.$

{\bf Theorem 12.} Let $P\in \mathcal{B(H)}^{Id}.$ If $P$ and $I-P$ have the operator matrices (2.14), respectively, then

(i) \ \ $P_{(2I-P-P^{\ast})^{+}}=P_{(P+P^{\ast})^{-}}+P_{N(P+P^*)}.$

(ii)\ \ $P_{N(P_{1}^{\ast})}=P_{N(Q_{1})}$ and  $P_{N(P_{1})}=P_{N(Q_{1}^{\ast})}.$

(iii) $P_{(2I-P-P^{\ast})^{-}}+P_{N(2I-P-P^{\ast})}=P_{(P+P^{\ast})^{+}}.$

{\bf Proof.} (i) Replacing $P$ by $I-P$ in equation (2.3), we have
$$min\{J: J(I-P)\geqslant0,\ J=J^{\ast}=J^{-1}\}=2P_{(2I-P-P^{\ast})^{+}}-I.$$
Also, Lemma 11 implies that $$min\{J: J(I-P)\geqslant0,\ J=J^{\ast}=J^{-1}\}=min\{J: P^{\ast}JP\leqslant J,\ J=J^{\ast}=J^{-1}\}.$$
Then
\begin{equation}P_{(2I-P-P^{\ast})^{+}}=P_{(P+P^{\ast})^{-}}+P_{N(P+P^*)}\end{equation} follows from (i) of Theorem 7.

(ii) Replacing $P$ by $I-P$ in (i) of Theorem 8, we get that
$$max\{J: J(I-P)\geqslant0,\ J=J^{\ast}=J^{-1}\}=2P_{(2I-P-P^{\ast})^{+}}-I+2P_{N(2I-P-P^{\ast})}.$$
In view of Lemma 11 and (ii) of Theorem 7, we have
\begin{equation}P_{(2I-P-P^{\ast})^{+}}+P_{N(2I-P-P^{\ast})}=P_{(P+P^{\ast})^{-}}+P_{N(P-P^{\ast})}.\end{equation}
 Then equations (2.27) and (2.28) imply that
\begin{equation}P_{N(2I-P-P^{\ast})}=P_{N(P-P^{\ast})}-P_{N(P+P^{\ast})},\end{equation}
so $$P_{N(Q_{1})}=P_{N(2I-P-P^{\ast})}=P_{N(P-P^{\ast})}-P_{N(P+P^{\ast})}=P_{N(P_{1}^{\ast})}$$ follows from Lemma 6.
In a similar way, we have \begin{equation}P_{N(P_{1})}=P_{N(P+P^{\ast})}=P_{N(P^{\ast}-P)}-P_{N(2I-P-P^{\ast})}=P_{N(Q_{1}^{\ast})}.\end{equation}

(iii) follows from (i).
\qquad$\square$

{\bf Lemma 13.} ([2, Theorem 2.3])For a $J$-projection $P$, there exists uniquely a $J$-positive projection $Q,$
a $J$-negative projection $R$ such that $$P=Q+R,\ \ \  QR=RQ=0\ \ \ \hbox {and} \ \ \ QR^{\ast}=R^{\ast}Q=0.$$

The following result is obtained in [2,7].
However, our method is more concrete than those of [2,7].

{\bf Corollary 14.} ([2, Theorem 2.12]) or [7, Theorem 2]) For a $J$-projection $P$, there exists uniquely a $J$-contractive projection $E_{1}$ and
a $J$-expansive projection $E_{2}$ such that $$P=E_{1}E_{2}=E_{2}E_{1}=E_{1}+E_{2}-I
\ \ \ \hbox {and} \ \ \ E_{1}E_{2}^{\ast}=E_{2}^{\ast}E_{1}=E_{1}+E_{2}^{\ast}-I.$$

{\bf Proof.} Suppose that $P$ has operator matrix form (1.1).
It follows from the fact $P=JP^{\ast}J$ that $J$ has the form (2.1) in Lemma 2.
With respect to the space decomposition
$\mathcal{H}=R(P)\oplus R(P)^{\perp},$ we define $E_{1}$ and $E_{2}$ as the following operator matrices
$$E_{1}:=\left(\begin{array}{cc}I_{1}&P_1\dfrac{I_{2}+J_2}{2}\\\\0&\dfrac{I_{2}-J_2}{2}\end{array}\right) \  \ \ \ \hbox {and}\ \ \ \
E_{2}:=\left(\begin{array}{cc}I_{1}&P_1\dfrac{I_{2}-J_2}{2}\\\\0&\dfrac{I_{2}+J_2}{2}\end{array}\right),$$
where $I_{1}$ and $I_{2}$ are identity operarors on the subspaces of $R(P)$ and $R(P)^{\perp}$, respectively.

Clearly,\ $E_{1}^{2}=E_{1}$ and\ $E_{2}^{2}=E_{2}.$
Then a direct calculation yields that $$P=E_{1}E_{2}=E_{2}E_{1}=E_{1}+E_{2}-I
\ \ \ \hbox {and} \ \ \ E_{1}E_{2}^{\ast}=E_{2}^{\ast}E_{1}=E_{1}+E_{2}^{\ast}-I,$$
as $P_{1}^{\ast}P_{1}J_{2}=J_{2}P_{1}^{\ast}P_{1}.$
It is easy to verify that $$J(I-E_{1})=0\oplus(I+P_1^{\ast}P_1)^{\frac{1}{2}}\dfrac{I_{2}+J_2}{2}\geqslant 0$$ and
$$J(I-E_{2})=0\oplus(I+P_1^{\ast}P_1)^{\frac{1}{2}}\dfrac{J_2-I_{2}}{2}\leqslant 0.$$ Thus by Lemma 11,
we conclude that $E_{1}$ is a $J$-contractive projection and $E_{2}$ is a $J$-expansive projection as desired.

To show the uniqueness, we note that $I-P$ is also a $J$-projection.
  If $$P=F_{1}F_{2}=F_{2}F_{1}=F_{1}+F_{2}-I
\ \ \ \hbox {and} \ \ \ F_{1}F_{2}^{\ast}=F_{2}^{\ast}F_{1}=F_{1}+F_{2}^{\ast}-I,$$
 where $F_1$ is a $J$-contractive projection and $F_{2}$
is a $J$-expansive projection.
 Setting $Q:=I-F_{1}$  and $R:=I-F_{2},$  we get that
$$I-P=Q+R,\ \ \  QR=RQ=0\ \ \ \hbox {and} \ \ \ QR^{\ast}=R^{\ast}Q=0.$$
Furthermore, Lemma 11 implies that $Q$ is a $J$-positive projection
and $R$ is a $J$-negative projection, so
 $Q=I-E_{1}$ and $R=I-E_{2}$ follow from the uniqueness of Lemma 13.
 Thus $F_1=E_1$ and $F_2=E_2.$  \qquad $\square$

\end{document}